\documentclass[12pt]{article}

\usepackage{amsmath}
\usepackage{lipsum}
\usepackage{tikz-cd}

\usepackage{amsthm}
\usepackage{mathrsfs}	
\usepackage{amssymb}
\usepackage{amscd}
\usepackage{setspace}	
\usepackage{tikz}		
\usepackage{tikz-cd}
\usetikzlibrary{matrix}	
\usepackage{graphicx}	
\graphicspath{ {pics/} }	
\usepackage{wrapfig}	
\usepackage{enumerate}	
\usepackage{bbm}		

\usepackage[utf8]{inputenc}	
\usepackage[english]{babel}	
\usepackage{blindtext}		
\usepackage[affil-it]{authblk}	

\usepackage{faktor}			
\usepackage{xparse}			

\usepackage{hyperref}		
\usepackage[lite]{amsrefs}

\newenvironment{acknowledgements}{%
  \begin{abstract}
}{%
  \end{abstract}
}

\newtheorem{theorem}{Theorem}[section] 

\newtheorem{definition}[theorem]{Definition}

\newtheorem{lemma}[theorem]{Lemma}

\newtheorem{remark}[theorem]{Remark}

\newtheorem{corollary}[theorem]{Corollary}

\newtheorem{notation}[theorem]{Notation}

\newcommand{\acts}{\curvearrowright}
\newcommand{\asdim}{\text{asdim}}
\newcommand{\diam}{\text{diam}}

\newcommand{\Z}{\mathbb{Z}}

\newcommand{\dad}{\text{DAD}}

\newtheorem*{theorem*}{Example}

\DeclareFontFamily{U}{mathb}{\hyphenchar\font45}
\DeclareFontShape{U}{mathb}{m}{n}{ <-6> matha5 <6-7> matha6 <7-8>
mathb7 <8-9> mathb8 <9-10> mathb9 <10-12> mathb10 <12-> mathb12 }{}
\DeclareSymbolFont{mathb}{U}{mathb}{m}{n}

\DeclareMathAccent{\abxring}{0}{mathb}{"38}

\DeclareFontFamily{U}{mathb}{\hyphenchar\font45}
\DeclareFontShape{U}{mathb}{m}{n}{ <-6> matha5 <6-7> matha6 <7-8>
mathb7 <8-9> mathb8 <9-10> mathb9 <10-12> mathb10 <12-> mathb12 }{}
\DeclareSymbolFont{mathb}{U}{mathb}{m}{n}

\begin{document}
\title{Topological Rigidity of the Dynamic Asymptotic Dimension}

\author{Samantha Pilgrim}

\maketitle

\paragraph{Abstract:}We show for a free action of a countable group $\Gamma$ on a finite-dimensional, compact metric space by homeomorphisms that the dynamic asymptotic dimension is either infinite or coincides with the asymptotic dimension of $\Gamma$.  

\section{Introduction}

The dynamic asymptotic dimension of a group action $\Gamma\acts X$ was first introduced in \cite{dasdimGWY}.  Although known to be related to dimension theories for group actions which take into account both the asymptotic dimension of $\Gamma$ and the covering dimension of $X$ \cite[Theorem 5.14]{kerr2017dimension}, there are many examples for which $\dad(\Gamma\acts X) = \asdim\Gamma$.  This was recently shown to always hold for finite dimensional actions on $0$-dimensional spaces (implied by a groupoid result in \cite{bonicke2023dynamic}); and for minimal actions by virtually cyclic groups \cite{bonicke2023nuclear} (the latter was also previously proven with the additional assumption of freeness in \cite{amini2020dynamic}).  However, a proof of this bound for general actions on higher-dimensional spaces has until now been elusive.  

In this note, we will show that $\dad(\Gamma\acts X)\in \{\asdim\Gamma, \infty\}$ for all free actions of countable groups by homeomorphisms on finite-dimensional, compact metric spaces.  The key ingredients of the proof are some topological rigidity, which we describe in $\ref{trick}$, and the inductive dimension.  More generally, we will also see how these new methods motivate a broader investigation of how small boundary-type properties can be applied to the $\dad$.  

Sharp bounds for the $\dad$ have been sought after for some time and are related to questions formulated by Willett and others (see \cite[Section 8]{warpedcones}).  The proof is also elementary and relatively natural, especially considering certain arguments in large-scale geometry.  While anticipated, this result is still somewhat surprising as the definition of $\dad$ involves open sets, and so one would expect it to be sensitive to the topology of $X$.  In fact, our main theorem will imply it is equivalent to define the $\dad$ using Borel or even arbitrary sets, at least for finite-dimensional actions covered by the theorem.  Additional connections with measure-theoretic properties were exposed in \cite{Conley2020BorelAD}.  

The dynamic asymptotic dimension and related theories have found relevance in a number of interesting problems, including the Farrell-Jones conjecture on manifold topology \cite[Section 4]{dasdimGWY}, calculations of $C^*$-algebra $K$-theory \cite{Guentner2016DynamicAD} \cite{bonicke2021dynamic}, and the classifiability of $C^*$-algebras \cite[section 8]{dasdimGWY} (see the introduction of \cite{https://doi.org/10.48550/arxiv.2201.03409}, for instance, for exposition on the classification program).  Sharper bounds for the dimension in particular have implications for the homology of group actions \cite[Theorem 3.36]{bonicke2021dynamic}, similar to how the homology of a manifold vanishes in degrees higher than its dimension.  

\section{Preliminaries}

We begin by fixing some notation.  Throughout, $\Gamma$ will be a countable, discrete group.  Such a group always admits a proper, right-invariant metric unique up to coarse equivalence \cite[Proposition 2.3.3]{notes_on_prop_A}; and we will fix, for all time, some such metric on $\Gamma$.  The space $X$ will always be compact and metrizable with some fixed metric, and $\Gamma\acts X$ will always be a free action by homeomorphisms.  We write $\gamma\cdot x$ for the action of $\gamma\in \Gamma$ on $x\in X$.  We denote by $B_x^r(X)$ the closed ball of radius $r$ in $X$ centered at $x$, and by $N_\epsilon(Y)$ the open $\epsilon$-neighborhood of a set $Y\subset X$.  

The asymptotic dimension of a metric space was originally introduced by Gromov in \cite{Gromov1991AsymptoticIO}.  There are many definitions, but we will use the one below.  

\definition{Let $G$ be a (in general non-compact) metric space and $A\subset G$ a subset.  An $r$-chain in $A$ is a finite sequence $x_1, \ldots, x_n$ of points in $A$ such that $d_G(x_i, x_{i+1})\leq r$ for $1\leq i\leq n-1$.  Two points $x, y\in A$ are in the same $r$-component of $A$ if they are connected by an $r$-chain in $A$.  We say $B\subset G$ is $R$-bounded if $\text{diam}(B)\leq R$.  A cover $\mathcal{U}$ of $G$ by $d+1$-sets such that the $r$-components of each $U\in \mathcal{U}$ are $R$-bounded is called a $(d, r, R)$-cover for $G$.  We say $\asdim G\leq d$ if for all $r>0$ there exists $R>0$ and a $(d, r, R)$-cover for $G$.  In this case, we call the function $D_G(r):= R$ a $d$-dimensional control function for $G$.  The asymptotic dimension of $G$ is then the least integer $d$ such that $\asdim G\leq d$ and is infinite if no such $d$ exists.  Since the groups $\Gamma$ we consider always carry a metric unique up to coarse equivalence, and asymptotic dimension is invariant under coarse equivalence, it is well-defined to write $\asdim\Gamma$ for the asymptotic dimension of $\Gamma$ with its fixed metric. }

\noindent \normalfont We will need to talk about ``partially defined actions" which arise in our setting from restricting an action $\Gamma\acts X$ to a subspace $Y\subset X$.  Since most of our notation and lemmas are not much harder to state for general partial dyanmical systems, we set things up in this level of generality.  

\definition{A topological partial action of $\Gamma$ on a topological space $X$ is a pair $(\{D_\gamma\}_{\gamma\in \Gamma}, \{\theta_\gamma\}_{\gamma\in \Gamma})$ consisting of a collection $\{D_\gamma\}_{\gamma\in \Gamma}$ of subsets of $X$, and a collection $\{\theta_\gamma\}_{\gamma\in \Gamma}$ of homeomorphisms, $\theta_\gamma: D_{\gamma}\to D_{\gamma^{-1}}$ such that  

\begin{description}
\item(i) $D_e = X$, and $\theta_e$ is the identity map.  

\item(ii) $\theta_\gamma\circ \theta_\lambda\subseteq \theta_{\gamma\lambda}$, for all $\gamma, \lambda \in \Gamma$.  
\end{description}

We do not assume the sets $D_\gamma$ are open, even though it is more customary to do so (more on this later).  Here, the composition $\theta_\gamma\circ \theta_\lambda$ denotes the map whose domain is the set of all $x\in X$ for which $\theta_\gamma(\theta_\lambda(x))$ is defined.  In other words, this is the set $\theta_\lambda^{-1}(D_\gamma) = \theta_{\lambda}^{-1}(D_{\lambda^{-1}}\cap D_\gamma)$.  The symbol ``$\subseteq$" means that the function on the right-hand-side is an extension of the function on the left-hand-side.  Notice also that $\theta_{\gamma^{-1}} = \theta_{\gamma}^{-1}$.  

A partial dynamical system is a quadruple $(X, \Gamma, \{D_\gamma\}_{\gamma\in \Gamma}, \{\theta_\gamma\}_{\gamma\in \Gamma})$.  An action $\Gamma\acts X$ is then a partial dynamical system with $D_\gamma = X$ for all $\gamma\in \Gamma$.  A partial dynamical system is free if $\theta_\gamma(x) = \theta_{\lambda}(x)$ iff $\gamma = \lambda$.  

}

\notation{The partial dynamical systems we will deal with arise from restricting actions, in which case we will write $\Gamma\acts A$ for a subspace $A\subset X$ to denote the partial dynamical system $(A, \Gamma, \{A\cap h_\gamma^{-1}(A)\}_{\gamma\in \Gamma}, \{h_\gamma\}_{\gamma\in \Gamma})$ where $h_\gamma$ comes from $\Gamma\acts X$.  Notice that if $A\subset X$ is closed, so are the sets $D_\gamma$, as they are obtained by intersecting closed sets.  The partial actions we will work with will always have closed domain sets $D_\gamma$.  }\label{restricted action notation}

\normalfont See \cite[Part I]{Exel2015PartialDS} for a more complete treatment of partial actions.  

\definition{If $(X, \Gamma, \{D_\gamma\}_{\gamma\in \Gamma}, \{\theta_\gamma\}_{\gamma\in \Gamma})$ is a free partial dynamical system on a compact space, $\dad_{free}(X, \Gamma, \{D_\gamma\}_{\gamma\in \Gamma}, \{\theta_\gamma\}_{\gamma\in \Gamma})$ is the smallest integer $d$ such that for all finite $F\subset \Gamma$ there is an open cover $\mathcal{U} = \{U_0,\ldots, U_d\}$ (that is, $U_i$ is open in $X$) of $X$ such that, for $0\leq j\leq d$ and $x\in U_j$, the set

$$Y_j^x = \Bigg\{\begin{array}{l|l} y\in X & \exists \text{ }  f_1, \ldots, f_k\in F\text{ such that } y = \theta_{f_k}\circ \cdots \circ \theta_{f_1}(x) \\ \text{ }& \text{ and } \theta_{f_{k_0}}\circ \cdots \circ \theta_{f_1}(x)\in U_j  \text{ (and exists) } \forall \text{ } 1\leq k_0\leq k \end{array}\Bigg\}$$
is uniformly finite (i.e. the bound on its cardinality is independent of $x$ and $j$).  For formal reasons, we further define that $\dad_{free}(X, \Gamma, \{D_\gamma\}_{\gamma\in \Gamma}, \{\theta_\gamma\}_{\gamma\in \Gamma}) = -1$ if $X = \emptyset$.As we will work only with free systems, we will now just write $\dad$.  This definition agrees with \cite[Definition 1.2]{dasdimGWY} for a free action $\Gamma\acts X$.  Note that if $\Gamma\acts Y\subset X$ is a partial dynamical system obtained by restricting an action $\Gamma\acts X$, then $\dad(\Gamma\acts Y)\leq \dad(\Gamma\acts X)$.  }\label{dad definition}  

%
%

\definition{Let $\mathcal{P}_{fs}(\Gamma)$ denote the collection of all finite subsets $F\subset \Gamma$ which contain the identity and are symmetric (meaning $F^{-1}=F$).  Let $(X, \Gamma, \{D_\gamma\}_{\gamma\in \Gamma}, \{\theta_\gamma\}_{\gamma\in \Gamma})$ be a partial dynamical system, $B\subset X$ a subset, and $F\in \mathcal{P}_{fs}(\Gamma)$.  An $F$-chain in $B$ is a finite sequence $x_0, \ldots, x_n$ of points in $B$ such that $x_{i+1} = \theta_{f_{i+1}}(x_i)$ with $f_{i+1}\in F$ ($0\leq i<n$).  Two points $x, y\in B$ are in the same $F$-component of $B$ if they are connected by an $F$-chain in $B$ (this is an equivalence relation since $F\in \mathcal{P}_{fs}(\Gamma)$).  If $B\subset V\subset X$, the $F$-component of $B$ in $V$ is the smallest subset of $V$ containing $B$ which contains every $F$-chain in $V$ beginning at a point in $B$.  If $A\subset X$ is another subset, we say $A$ and $B$ are $F$-separated if their $F$-components in $A\cup B$ are disjoint.  

If $S\in \mathcal{P}_{fs}(\Gamma)$, we say $Y\subset X$ is $S$-bounded if it is a subset of $\{\theta_s(x) | s\in S \text{ and } x\in D_s\}$ for some $x\in X$.  A cover $\mathcal{U}$ of $X$ by $d+1$ (open/closed) sets whose $F$-components are $S$ bounded is called a (open/closed) $(d, F, S)$-cover for $(X, \Gamma, \{D_\gamma\}_{\gamma\in \Gamma}, \{\theta_\gamma\}_{\gamma\in \Gamma})$ or simply for $X$ if unambiguous.  Notice that $S$ can always be taken to be $F^r$ for some $r>0$ so that $S\in \mathcal{P}_{fs}(\Gamma)$.  We will write $S_i^x$ for the subset of $S$ such that $\{\theta_{s}(x) | s\in S_i^x\}$ is the $F$-component of $x$ in $U_i\in \mathcal{U}$.    A $d$-dimensional control function for $(X, \Gamma, \{D_\gamma\}_{\gamma\in \Gamma}, \{\theta_\gamma\}_{\gamma\in \Gamma})$ is a function $\mathcal{D}_{X}: \mathcal{P}_{fs}(\Gamma)\to \mathcal{P}_{fs}(\Gamma)$ such that for every $F\in \mathcal{P}_{fs}(\Gamma)$, there is a $(d, F, \mathcal{D}_X(F))$-cover for $(X, \Gamma, \{D_\gamma\}_{\gamma\in \Gamma}, \{\theta_\gamma\}_{\gamma\in \Gamma})$.  Note that $\dad_{free}(X, \Gamma, \{D_\gamma\}_{\gamma\in \Gamma}, \{\theta_\gamma\}_{\gamma\in \Gamma})\leq d$ iff there is a $d$-dimensional control function for $(X, \Gamma, \{D_\gamma\}_{\gamma\in \Gamma}, \{\theta_\gamma\}_{\gamma\in \Gamma})$.  } \label{chain definition}

\normalfont \noindent An important ingredient in \cite[Proposition 3.1]{boxspacesDT} is a finite union theorem for the asymptotic dimension.  We will prove an analogous result for the dynamic asymptotic dimension, inspired by \cite[Lemma 3.5]{Brodskiy2006AHT}, and similar to \cite[Theorem 2.6]{bonicke2023dynamic}.  

\lemma{(Finite union lemma) Let $(X, \Gamma, \{D_\gamma\}_{\gamma\in \Gamma}, \{\theta_\gamma\}_{\gamma\in \Gamma})$ be a partial dynamical system.  Let $A, B\subset X$ and $F\in \mathcal{P}_{fs}(\Gamma)$.  Assume that the $F^{r_A}$-components of $A$ are $F^{R_A}$-bounded, that the $F^{r_B}$-components of $B$ are $F^{R_B}$-bounded, and that $2R_B + 2r_B<r_A$.  Then the $F^{r_B}$-components of $A\cup B$ are $F^{r_A + R_A}$-bounded.  }\label{finite union lemma}
\begin{proof}

Let $x_0, \ldots, x_n$ be an $F^{r_B}$-chain in $A\cup B$.  Suppose $x_j$ and $x_k$ are two consecutive points contained only in $A$.  Then $x_{j+1}, \ldots, x_{k-1}$ is an $F^{r_B}$-chain in $B$ and therefore is $F^{R_B}$-bounded.  

We therefore have that $x_{j+1}$ and $x_{k-1}$ are in the same $F^{R_B}$-component, so $x_i$ and $x_j$ are in the same $F^{2R_B + 2r_B}$-component of $A$.  Since $2R_B + 2r_B < r_A$, the points in the original chain which are contained only in $A$ therefore form an $F^{r_A}$-chain in $A$ which is then $F^{R_A}$-bounded.  The original chain is therefore $F^{2R_B + 2r_B + R_A}$-bounded.  So the $F^{r_B}$-components of $A\cup B$ are $F^{2R_B + 2r_B + R_A}$-bounded, which implies the result since $2R_B + 2r_B< r_A$.  \end{proof}

\lemma{Suppose $\asdim\Gamma\leq d$ with control function $D_\Gamma$.  Then there is a $d$-dimensional control function $\mathcal{D}_\Gamma$ for $\Gamma\acts \Gamma$ given by $F\mapsto B_e^{D_\Gamma(\text{diam}(F))}(\Gamma)$ ($F\in \mathcal{P}_{fs}(\Gamma)$).  }\label{basic lemma}
\qed

\normalfont \noindent The following lemma is of general interest.  We will use the second part just once in our main proof, but will use the first part multiple times.  

%
%

\lemma{Let $\Gamma\acts X$ be a free action of a countable group by homeomorphisms on a compact metric space.  Suppose $U\subset X$ is closed and that the $F$-components of $U$ are $S$-bounded.  For $x\in U$, let $S^x\subset \Gamma$ be the unique set such that $S^x\cdot x$ is the $F$-component of $x$ in $U$.  Then

\begin{description}
\item(i) $U$ has an open neighborhood in $X$ whose $F$-components are $S$-bounded.  

\item(ii) About each $x\in U$, there is an open neighborhood $U^x$ in $X$ such that if $s\in S^x$, $f\in F$, but $fs\notin S^x$, then $d(fs\cdot U^x, U)>0$.  

\end{description}} \label{trick}

\begin{proof}

We can assume $S\in \mathcal{P}_{fs}(\Gamma)$.  To prove $(i)$, suppose not, that is, suppose for every $n$, the $F$ components of $N_{1/n}(U)$ are not $S$-bounded.  So for every $n$, we have an $F$-chain $x^{(n)}_0, \ldots, x^{(n)}_{K_n}$ in $N_{1/n}(U)$ with $x^{(n)}_0, \ldots, x^{(n)}_{K_n-1}\in S\cdot x_0$ but $x_{K_n} = \gamma_n\cdot x^{(n)}_0$ for $\gamma_n = f_ns_n$ with $s_n\in S$, $f_n\in F$, and $\gamma_n\notin S$.  Since there are finitely-many $\gamma\in \Gamma$ with the property that $\gamma = fs$ for $s\in S$ and $f\in F$, we must have $\gamma_n = \lambda$ for infinitely-many $n$.  Furthermore, there are finitely-many sequences $f_1, \ldots, f_{K_n}$ ($f_k\in F$) such that $f_k\cdots f_1\in S$ for all $1\leq k<K_n$ but $f_{K_n}\cdots f_1 = \lambda\notin S$, and so for infinitely-many $n$ we have an $F$-chain $x^{(n)}_0, \ldots, x^{(n)}_{K}$ in $N_{1/n}(U)$ and $f_1, \ldots, f_K\in F$ such that $f_{k+1}\cdot x^{(n)}_{k} = x^{(n)}_{k+1}$ for $0\leq k<K$, and also $x_K = \lambda\cdot x_0$ with $\lambda\notin S$.  Since $X$ is compact, we have a subsequence $x_0^{(n_j)}$ converging to a point $y_0\in X$, and so $f_k\ldots f_1\cdot x_0^{(n_j)} = x_k^{(n_j)}$ converges to $f_k\ldots f_1\cdot y_0$ for all $1\leq k\leq K$ by continuity.  Since $x_k^{(n_j)}\in N_{1/n_j}(U)$ and $U$ is closed in $X$, $x_k^{(n_j)}$ converges to a point in $U$.  Hence, we have an $F$-chain in $U$ which is not $S$-bounded.  



To prove $(ii)$, let $Y^x$ be the $F$-component of $x\in U$.  We have

$$\delta_x := \min_{y\in Y^x} \min_{\{f\in F | f\cdot y\notin U\}} d_X(U, f\cdot y)>0.$$
Let $0<\eta_x$ such that $d(z, x)<\eta_x$ implies $d(\gamma\cdot z, \gamma\cdot x)<\delta_x/2$ for all $\gamma\in FS^x$.  Then putting $U^x = \text{Int}(B_{x}^{\eta_x}(X))$ proves $(ii)$.  \end{proof}

\lemma{Suppose $\Gamma\acts X$ is free and $\dad(\Gamma\acts X)\leq d$.  Then $\asdim\Gamma\leq d$.  }\label{easy inequality}

\begin{proof}
Since $\Gamma\acts X$ is free, we can identify $\Gamma$ with the orbit of any $x\in X$.  Then a $(d, B_e^R(\Gamma), S)$-cover for $\Gamma\acts X$ gives rise to a $(d, R, \diam(S))$-cover of $\Gamma$.  In fact, this works if the cover of $X$ we start with uses Borel or even arbitrary sets. \end{proof}

\normalfont \noindent The reverse inequality is more difficult.  



\section{Sharp Bounds for the DAD}

\normalfont
We will now give a short exposition of our methods before the main argument.  When $X$ is a Cantor set, the proof that $\dad(\Gamma\acts X)\leq \asdim\Gamma$ whenever $\dad(\Gamma\acts X)<\infty$ is fairly natural given the proof of \cite[Proposition 3.1]{boxspacesDT}.  Being somewhat imprecise, the idea is to use the assumption that $\dad(\Gamma\acts X)<\infty$ to break $X$ into finitely-many pieces each with finite components at an appropriate scale (i.e. for a given $F\in \mathcal{P}_{fs}(\Gamma)$).  We can then cover each piece with disjoint clopen towers (see \cite{kerr2017dimension} for a definition), which are separated and each have dimension at most $\asdim\Gamma$ at the given scale because the action restricted to a tower will resemble the structure of $\Gamma$ at that scale.  Covers for these towers can be combined at no cost because they are separated, and the finite union lemma allows us to put the pieces back together.  In higher dimensions, towers cannot be both open and disjoint unless they do not cover all of $X$.  However, one can arrange things so that the remaining part of $X$ is contained in the boundaries of some sufficiently nice sets, and so is one dimension lower.  This allows a pleasing argument using the inductive dimension.  

Aside from the many technicalities involved in making this idea precise, the main obstacle is that we can't commute taking boundaries with translation by a partial action in the way one might want, as the boundaries of the domain sets can come into play.  This appears to be related to the requirement in \cite[Corollary 2.15]{bonicke2023dynamic} that the domain sets be clopen, as this is equivalent to requiring that their boundaries be $(-1)$-dimensional (see the definition below).  In the setting of groupoids, it is also related to how restricting to a subspace of the unit space does not in general preserve the \'{e}tale condition.  We are grateful to Jeremy Hume for pointing out this fact about groupoids.  

The methods employed in the proof of our main theorem motivate a broader investigation of small boundary-type properties and their relation to the dynamic asymptotic dimension.  Indeed, one can probably give a similar and possibly  simpler proof of \ref{main corollary} using the topological small boundary property in place of the inductive dimension, though it requires an application of \cite[Theorem 3.8]{10.1112/plms/pdu065}, and so requires the space $X$ to be finite dimensional, which we don't assume in \ref{main theorem}.  It is worth remarking further that any compact metrizable space is an increasing union of finite-dimensional, closed subspaces, so we may be close to discharging the assumption that $X$ is finite-dimensional.

\noindent The following definition can be found in \cite[Chapter 7, Section 1]{Eng89}, but the concept of inductive dimension is older.

\definition{We define the (large) inductive dimension as follows.  Start by defining $\dim_{Ind}(X) = -1$ iff $X=\emptyset$.  Then define that $X$ has $\dim_{Ind}(X)=d$ for $d\geq 0$ if $d$ is the smallest integer such that for every closed subset $Y\subset X$ and open neighborhood $U$ of $Y$, there is an open set $U'$ with $Y\subset U'\subset U$ and $\dim_{Ind}(\partial U')\leq d-1$.  }

\remark{If $\dim_{Ind}(X) \leq d$, $\dim_{Ind}(Y)\leq d$ for any subspace $Y\subset X$; a finite union of spaces of dimension at most $d$ has dimension at most $d$; the image of a $d$-dimensional space under a homeomorphism is $d$-dimensional; and a countable union of compact spaces with dimension at most $d$ has dimension at most $d$.  The inductive dimension agrees with the Lebesgue covering dimension when $X$ is separable metric space \cite[Theorem 7.3.3]{Eng89}, so in particular when $X$ is a compact metric space.  Separability is equivalent to second countability for metric spaces, and so always passes to subspaces.}

\theorem{Let $X$ be a compact metric space and $Y\subset X$ a closed subspace with $\dim_{Ind}(Y)\leq m$.  Suppose $\Gamma\acts X$ is a free action by homeomorphisms and that $\dad(\Gamma\acts X)\leq d$ for some $d$.  Then $\dad(\Gamma\acts Y)\leq \asdim\Gamma$ where $\Gamma\acts Y$ is the partial dynamical system obtained by restricting $\Gamma\acts X$ to $Y$ (see \ref{restricted action notation}).  }\label{main theorem}

\begin{proof}
We will proceed by induction on the dimension of $Y$.  Assume the theorem holds when $Y$ has  dimension $\leq m-1$.  Let $\mathcal{D}_\Gamma$ be the $\asdim\Gamma$-dimensional control function for $\Gamma\acts \Gamma$ coming from \ref{basic lemma}.  We can assume $\mathcal{D}_\Gamma(F)\supset F$ for all $F\in \mathcal{P}_{fs}(\Gamma)$.  Fix $F_0\in \mathcal{P}_{fs}(\Gamma)$, let $G_0 = \mathcal{D}_\Gamma(F_0)$, and define $G_i$ and $F_i$ for $i = 1, \ldots, d$ by the inductive formulae $F_{i} = G_{i-1}^4$ and $G_{i} = \mathcal{D}_\Gamma(F_{i})^2$.  Except when above the letters $G$ or $T$, superscripts will be indices, not exponents.  Note each $F_i$ is symmetric and that $F_{i+1}\supset F_i$.  

We will continue using $\gamma\cdot x$ only for the action $\Gamma\acts X$.  Start with an open $(d, F_d, S)$-cover $\mathcal{U} = \{U_0, \ldots, U_d\}$ for $\Gamma\acts X$.  This cover has some lebesgue number $\lambda$, so if $U'_i$ is the closure of the $\lambda/2$-interior of $U_i$ and $U''_i$ is the closure of the $4\lambda/5$-interior of $U_i$, then $\mathcal{U}' = \{U_0', \ldots, U_d'\}$ and $\mathcal{U}'' = \{U_0'', \ldots, U_d''\}$ are closed $(d, F_d, S)$-covers for $\Gamma\acts X$.  


Put $\widehat{Y} = \Gamma\cdot Y$ (note that this set is not necessarily closed in $X$).  Since $Y$ is closed, hence compact, $\dim_{Ind}(\widehat{Y})\leq m$.  If we let $V_i = U_i''\cap \widehat{Y}$, then $\mathcal{V} = \{V_0, \ldots, V_d\}$ is a closed $(d, F_d, S)$-cover for $\widehat{Y}$.  Using the large inductive dimension, we can find open sets $V_i'$ in $\widehat{Y}$ such that $V_i\subset V'_i\subset N_{\lambda/5}(V_i)\cap \widehat{Y}\subset U_i'$ and $\dim_{Ind}(\partial_{\widehat{Y}} V'_i)\leq m-1$; and if we let $\widehat{V}_i$ be the closure of $V_i'$ in $\widehat{Y}$, then still $\dim_{Ind}(\partial_{\widehat{Y}} \widehat{V}_i)\leq m-1$ (since $\partial_{\widehat{Y}} \widehat{V}_i\subset \partial_{\widehat{Y}} V'_i$), and $\widehat{\mathcal{V}} = \{\widehat{V}_0, \ldots, \widehat{V}_d\}$ is a closed $(d, F_d, S)$-cover of $\widehat{Y}$.  Note that $\widehat{V}_i\subset U_i'$.  

For $x\in U'_i$, let $S_i^x\subset \Gamma$ be the unique set such that $S_i^x\cdot x$ is the $F_d$-component of $x$ in $U'_i$.  Using $\ref{trick}$(ii), we can find about each $x\in U_i'$ a neighborhood $U_i^x$ in $X$ such that if $s\in S_i^x$, $f\in F_d$, but $fs\notin S_i^x$, then $d(fs\cdot U_i^x, U_i')>0$.  Now find an open (in $\widehat{Y}$) neighborhood $W_i^y\subset U_i^y$ about each $y$ in $\widehat{V}_i$ such that $\dim_{Ind}(\partial_{\widehat{Y}}W_i^y)\leq m-1$, and such that $d(s\cdot W_i^y, s'\cdot W_i^y)>0$ for all $s, s'\in S$ with $s\neq s'$ (this is possible since the action is free).  

Notice for $A\subset \widehat{Y}$ that $\partial_{\widehat{Y}}(\gamma\cdot A) = \gamma\cdot \partial_{\widehat{Y}} A$ (this relies on the fact that $\widehat{Y}$ is $\Gamma$-invariant).  Thus, $\partial_{\widehat{Y}} (\gamma\cdot (A\cap \widehat{V}_i)\cap \widehat{V}_i)\subset \gamma\cdot \partial_{\widehat{Y}} A\cup \gamma\cdot \partial_{\widehat{Y}} \widehat{V}_i\cup \partial_{\widehat{Y}} \widehat{V}_i$.  Let $C_i^y$ denote the $F_d$-component of $W_i^y\cap \widehat{V}_i$ in $\widehat{V}_i$.  Then $C_i^y$ can be described by starting with $W_i^y\cap \widehat{V}_i$, translating by an element of $F_d$, intersecting with $\widehat{V}_i$, translating by another element of $F_d$, intersecting again, and so on.  Since $W_i^y\subset U_i^y$, $\partial_{\widehat{Y}} C_i^y\subset S\cdot \partial_{\widehat{Y}} W_i^y\cup S\cdot \partial_{\widehat{Y}} \widehat{V}_i$ (recall that $S$ contains the identity), and so $\dim_{Ind}(\partial_{\widehat{Y}} C_i^y)\leq m-1$.  Since $C_i^y$ contains an open neighborhood about $y$ in $\widehat{Y}\cap \widehat{V}_i$, compactness of $Y\cap \widehat{V}_i$ implies we have a finite collection $\{C_i^{y_k}\}_{k=1}^{K_i}$ which covers $Y\cap \widehat{V}_i$.  To simplify notation, we will now simply write $C_i^k$ for $C_i^{y_k}$.  

Put $D_i^1 = C_i^1$ and $D_i^k = C_i^k\setminus \cup_{j<k} D_i^j = C_i^k\setminus \cup_{j<k} C_i^j$ for $k = 2, \ldots, K_i$.  By construction, $D_i^1$ is equal to its own $F_d$-component in $\widehat{V}_i$ and is disjoint from $D_i^2\subset \widehat{V}_i$.  Using that $F_d$ is symmetric, we therefore have that $D_i^1$ and $D_i^2$ are $F_d$-separated.  But then $D_i^2$ is equal to its own $F_d$-component in $\widehat{V}_i$ and is disjoint from $D_i^3\subset \widehat{V}_i^Y$, so by repeating the latter argument finitely-many times, we see that $D_i^k$ and $D_i^{k'}$ are $F_d$-separated for $k\neq k'$.  Notice also that $\cup_k D_i^k\supset \widehat{V}_i\cap Y$.

Let $Z = (\cup_i \cup_k \partial_{\widehat{Y}} C_i^k)\cap Y$.  Then $Z$ is closed in $Y$ and therefore also closed in $X$, and $\dim_{Ind}(Z)\leq m-1$.  Notice that $(\partial_{\widehat{Y}}D_i^k)\cap Y\subset Z$ for all $i$ and $k$ since in general, $\partial(A\setminus B)\subset \partial A\cup\partial B$ and $\partial(A\cup B)\subset \partial A\cup \partial B$.  By the inductive hypothesis and the assumption that $\Gamma\acts X$ is finite-dimensional, we can find an $(\asdim\Gamma, G_d^4, T)$-cover for $Z$ by sets open in the subspace topology on $Z$ (and we can assume without loss of generality that $G_d^4\subset T$).  Since $Z$ is closed in $Y$ and $Y$ is closed in $X$, $Z$ is closed in $X$ and therefore compact, and so this cover has some lebesgue number $\mu$, and so we can replace the sets in the cover by the closures of their $\mu/2$-interiors to get a closed $(\asdim\Gamma, G_d^4, T)$-cover for $Z$.  Then by part (i) of \ref{trick}, we can again replace this cover by a cover using sets which are open in $X$ and so cover an open  neighborhood $N$ of $Z$ in $X$.

For each $i$ and $k$, $D_i^k\setminus \partial_{\widehat{Y}}D_i^k$ is open in $\widehat{Y}$, and so $(D_i^k\setminus \partial_{\widehat{Y}} D_i^k)\cap Y$ is open in $Y$.  Moreover, since we chose $W_i^y$ so that its translates by different elements of $S$ are separated, $D_i^k\setminus \partial_{\widehat{Y}}D_i^k$ decomposes into a separated collection $\{E_i^{k, s} = (D_i^k\setminus \partial_{\widehat{Y}}D_i^k)\cap s\cdot W_i^{y_k} | s\in S'\}$ of non-empty open sets for some $S'\subset S_i^{y_k}\subset S$ with the property that $f\cdot E_i^{k, s}\subset E_i^{k, fs}\cup Z\cup \widehat{V}_i^c$ for $f\in F_d$, so long as we interpret $E_i^{k, fs}$ as being empty whenever $fs\notin S'$ (this is like being a tower in a very weak sense).  These properties imply an $(\asdim\Gamma, F_i, \mathcal{D}_\Gamma(F_i))$-cover for $\Gamma\acts \Gamma$ gives rise, using the bijection $\{E_i^{k, s} : s\in S'\}\leftrightarrow S'\subset \Gamma$, to an open $(\asdim\Gamma, F_i, \mathcal{D}_\Gamma(F_i))$-cover for $D_i^k\setminus \partial_{\widehat{Y}}D_i^k$ for each $k$.  Since $D_i^k$ and $D_i^{k'}$ are $F_d$-separated for $k\neq k'$, we can put these covers together at no cost to form an open $(\asdim\Gamma, F_i, \mathcal{D}_\Gamma(F_i))$-cover for $D_i := \cup_k (D_i^k\setminus \partial_{\widehat{Y}}D_i^k)$.  Notice that $D_i\supset (\widehat{V}_i\cap Y)\setminus Z$.  By our construction of $F_i$ and $G_i$, and $d(\asdim\Gamma + 1)$ applications of \ref{finite union lemma}, we can put the covers of the $D_i$ together to get an open (in $\widehat{Y}$) $(\asdim\Gamma, F_0, G_d)$-cover for $D := \cup_i D_i$.  

We have that $D\cup N \supset Y$.  Using the open (in $X$) $(\asdim\Gamma, G_d^4, T)$-cover for $N$ we found earlier and the open (in $\widehat{Y}$) cover for $D$ we found above, we can apply \ref{finite union lemma} $\asdim\Gamma + 1$ more times to obtain an open (in $Y$) $(\asdim\Gamma, F_0, T^2)$-cover for $Y$ so that we have shown $\dad(\Gamma\acts Y)\leq \asdim\Gamma$.  In the case where $m=0$, there is nothing more to do as $\partial_{\widehat{Y}} C_i^k = \emptyset$, and so $Z = \emptyset$.  \end{proof}

\corollary{If $X$ is a compact metric space with finite covering dimension and $\Gamma\acts X$ is a free action of a countable group by homeomorphisms, then $\dad(\Gamma\acts X) \in \{\asdim\Gamma, \infty\}$.} \label{main corollary}
\begin{proof} As remarked earlier, when $X$ is a compact metric space, the inductive dimension coincides with the Lebesgue covering dimension.  The corollary then follows immediately from the previous theorem and \ref{easy inequality}.  \end{proof}


\corollary{If $X$ is a compact metric space with finite covering dimension and $\Gamma\acts X$ is a free action of a countable, virtually-nilpotent group by homeomorphisms, then $\dad(\Gamma\acts X) = \asdim\Gamma$.  In particular, $\dad(\Gamma\acts X) = d$ if $\Gamma$ is virtually isomorphic to $\Z^d$.  } 
\begin{proof}
We only need to know that $\dad(\Gamma\acts X)<\infty$, which follows from \cite[Theorem 8.5]{warpedcones}.  \end{proof}

\noindent \normalfont We can now give a complete description of the dynamic asymptotic dimension in terms of intrinsic properties of the acting group for a reasonably large class of examples.  

\corollary{If $X$ is a compact Lie group and $\Gamma<X$ is a countable subgroup, then the translation action $\Gamma\acts X$ is infinite-dimensional if $\Gamma$ is not amenable and has $\dad(\Gamma\acts X) = \asdim\Gamma$ if $\Gamma$ is amenable.}
\begin{proof}
Since the asymptotic dimension of an increasing union of groups is the supremum over the dimensions of those groups, we can assume $\Gamma$ is finitely-generated.  Since $\Gamma\acts X$ preserves the Haar measure, $\dad(\Gamma\acts X) = \infty$ whenever $\Gamma$ is not amenable by \cite[Corollary 8.27]{dasdimGWY}.  On the other hand, a combination of representation theory, Lie's theorem, and Engel's theorem shows a finitely-generated, amenable subgroup of a compact Lie group is virtually abelian (see \cite{stacknew} for a proof), so if $\Gamma$ is amenable, $\dad(\Gamma\acts X) = \asdim\Gamma$.  \end{proof}

\begin{acknowledgements}
\noindent \raggedright This project has received funding from the European Research Council (ERC) under the European Union's Horizon 2020 research and innovation programme.  
\end{acknowledgements}



%

\bibliography{mybibliography2.bib}
\bibliographystyle{plain}

\end{document}